\documentclass{article}
\usepackage[T1]{fontenc}
\usepackage{lmodern}
\usepackage[a4paper]{geometry}
\usepackage{graphicx}
\usepackage{color}
\usepackage{amsmath}
\usepackage{amsfonts}
\usepackage{bbm}

\def\diffd{\mathrm{d}}

\DeclareMathOperator{\E}{\mathbb{E}}
\let\P\relax\DeclareMathOperator{\P}{\mathbb{P}}
\DeclareMathOperator{\one}{\mathbbm{1}}

\DeclareMathOperator{\arcsinh}{sinh^{-1}}
\DeclareMathOperator{\cotanh}{cotanh}

\def\all{a}
\def\majo{M}
\def\mino{m}
\def\sigO{{(0,0,\ldots,0)}}
\def\sigI{{(1,1,\ldots,1)}}
\def\sigend{\sigma_\text{fittest}}

\newtheorem{theorem}{Theorem}
\newtheorem{conjecture}{Conjecture}
\newtheorem{lemma}{Lemma}

\title{Accessibility percolation with backsteps}
\author{Julien Berestycki,
\'Eric Brunet, Zhan Shi\footnote{%
emails: \texttt{Julien.Berestycki@upmc.fr},
\texttt{Zhan.Shi@upmc.fr},
\texttt{EricBrunet@lps.ens.fr}.\newline
J.B. and Z.S.: Sorbonne Universit\'es, UPMC Univ Paris 06, CNRS, UMR 7599,
LPMA, F-75005, Paris France\newline
\'E.B.: Sorbonne Universit\'es, UPMC Univ Paris 06, CNRS, UMR 8550, LPS-ENS,
F-75005, Paris France
}}

\begin{document}

\maketitle

\begin{abstract}
Consider a graph in which each site is endowed with  a value 
called \emph{fitness}. A path in the graph is said to be
``open'' or ``accessible'' if the fitness values along that path is strictly increasing. We say
that there is accessibility percolation between two sites when such a path
between them exists. Motivated by the so called House-of-Cards model from
evolutionary biology, we consider this question on the $L$-hypercube
$\{0,1\}^L$ where the fitness values are independent random variables. We show
that, in the large $L$ limit, the probability that an accessible path
exists from an arbitrary starting point to the (random) fittest site is no
more than $x^*_{1/2}= 1-\frac12\sinh^{-1}(2) =0.27818\ldots$ and we
conjecture that this
probability does converge to $x^*_{1/2}$. More precisely, there is a phase
transition on the value of the fitness~$x$ of the starting site: assuming that
the fitnesses are uniform in $[0,1]$, we show that, in the large $L$ limit,
there is almost surely no path to the fittest site  if $x>x^*_{1/2}$ and we
conjecture that there are almost surely many paths if $x<x^*_{1/2}$.
 If one conditions on the fittest
site to be on the opposite corner of the starting site rather than being
randomly chosen, the picture remains the same but with the critical point
being now $x^*_1= 1-\sinh^{-1}(1)= 0.11863\ldots$.
Along the way, we obtain a large~$L$ estimation for the number of
self-avoiding paths joining two opposite corners of the $L$-hypercube.

\medskip

\end{abstract}

\section{Introduction}
\subsection{Definition of the model}

We consider the following mathematical model inspired by evolutionary biology:
\begin{enumerate}
\item
The genome of an organism is made of $L$ sites which can each be in two
states (or alleles): 0 (the wild state) or 1 (the mutant state). There are
therefore $2^L$ possible genomes, which are coded as an $L$-bit binary word
or as a corner of the $L$-hypercube \cite{kauffman-levin}.
\item During reproduction (supposed asexual and without recombination), we
assume that the only mutations that can occur consist in changing the state
at one single site, either from the wild to the mutant state or from the
mutant to the wild state. With our representation, a mutation is flipping
one single bit in the $L$-bit word or traveling along one edge of the
$L$-hypercube \cite{klg,kauffman-levin}.
\item We assume that we are in a regime with a low mutation rate, high
selection and a population which is not too large. In this regime, when
a mutation occurs, the new genome either fixates (\textit{i.e.}\@ it
invades the whole population and becomes the new resident type) if its
fitness value is better than the value of the resident population, or is
eliminated if it is lower. This happens (in the regime we assume)
fast enough that a new mutation has no time to appear before the population
is homogeneous again.
\end{enumerate}
In this model, the evolutionary history of the population as a whole can be
described as a path along the edges of the $L$-hypercube, with the
constraint that the fitness value must increase at each step. We call such paths
``open'' or ``selectively accessible'' \cite{franke-kloezer-devisser-krug,
weinreich-watson-chao, weinreich-delaney-depristo-hartl}. We
emphasize that we allow paths of arbitrary length, where bits can flip from
1 to 0 as well as from 0 to 1. The question we
wish to address is the following: assuming that the population is initially
in the state $\sigO$, is there an evolutionary  path allowing it to evolve
to the fittest site available?

To answer this question we need a model for the fitness values of each site.
As a first approach, we consider the House-of-Cards \cite{kingman} model
(which is equivalent \cite{altenberg} to the $NK$ model
\cite{kauffman-levin} with
$K=N-1$) where the fitness values of the $2^L$ sites are independent
random numbers. For the purpose of discussing the existence of
open paths, the actual fitness values
 of each site are not relevant; the only useful
information are how the fitness values are ordered. This means that the answer
to our question does not depend on the chosen distribution of the
fitness values, and that we can safely choose the most convenient distribution:
\begin{enumerate}\setcounter{enumi}3
\item The fitness value of the fittest site is 1 and that the fitness
values of the other $2^L-1$ sites are independent random numbers chosen
uniformly between 0 and~1.
\end{enumerate}

As we explained, this is equivalent to the House-of-Cards model if,
furthermore, the fittest site is chosen
uniformly at random  amongst the $2^L$ sites of the hypercube. In this
paper, we first consider the case where the fittest site is
deterministically chosen to be $\sigI$, then the case where it is a fixed
arbitrary site $\sigend$ and, finally, the case where the fittest
site is random.

\subsection{Notations}
\begin{itemize}
\item Sites are coded as $L$ bit binary word. The initial state of the
population is $\sigO$.
\item $H$ (as in  ``Hamming distance'') is the number of bits set to 1 in
the fittest site $\sigend$.
\item The fitness of the starting site $\sigO$ is noted $x$.
\item $\Theta$ is the number of open (selectively accessible) paths from
$\sigO$ to the fittest site. To compare with previous results, we also
write as $\tilde\Theta$ the number of open paths of minimal length to the
fittest site.
\item The probability of an event is written $\P$ and its expectation $\E$. We
often need to condition on the value $x$ of the fitness in the starting
site; when we do we write the conditional probability and expectation
as $\P^x$ and $\E^x$. The values of $H$ and $L$ are implicit in the notation.
\end{itemize}

\subsection{Results of previous works}

Similar models have been studied by several groups in the past few years,
either directly on the hypercube as above \cite{us, hegarty-martinsson,
DePristo.2007, schmiegelt-krug} or in the geometrically
simpler setting of a tree \cite{us, roberts-zhao, nowak-krug}. 

Except for \cite{DePristo.2007}, all the previous studies focused only on
the number $\tilde\Theta$ of open paths going from the
starting position $\sigO$ to the opposite corner $\sigI$
 with minimal length, meaning that a mutation can only flip a bit
from 0 to 1 and not the other way around.
In this setting, one only needs to consider $H=L$ as direct paths to
a fittest site at Hamming distance $H$ cannot leave anyway the
$H$-hypercube .

There are $L!$ minimal length paths (open or not) connecting the starting
site $\sigO$ to the opposite corner $\sigI$, each of these minimal length 
paths go through $L$ random fitnesses between 0 and 1 (including the starting
site, but excluding the end site which is assumed to have fitness 1) and
the probability that a given minimal length path is open is the probability
that these $L$ random numbers are in order, which is $1/L!$. Therefore
\begin{equation}
\E(\tilde\Theta)=1.
\end{equation}
This expectation is however misleading, as the typical number of open
minimal length paths is not 1. Indeed, if one conditions on the fitness $x$
of the starting site, the probability that a given minimal length path is
open is $(1-x)^{L-1}/(L-1)!$ because the path meets $L-1$ random values
(excluding both the starting and end sites) and these values must be
all between $x$ and 1 and in ascending order. Therefore
\begin{equation}
\E^x(\tilde\Theta)=L(1-x)^{L-1}. \label{expminimal}
\end{equation}
This conditional expectation is a decreasing function of $x$ which is equal
to 1 for $x=x_c(L)$ with $x_c(L):=1-\exp[-\frac{\ln L}{L-1}]\sim\frac{\ln
L}L$ for large~$L$. This implies that
\begin{equation}
\P(\tilde\Theta\ge1)\le x_c(L)+\int_{x_c(L)}^1
\E^x(\tilde\Theta)=x_c(L)+(1-x_c(L))^L\sim \frac{\ln L}L
\qquad\text{ as }L\to\infty.
\end{equation}

By a clever second moment argument, Hegarty and Martinsson \cite{hegarty-martinsson}
proved that the above bound is tight:
\begin{equation}
\P(\tilde\Theta\ge1)\sim \frac{\ln L}L
\qquad\text{ as }L\to\infty.
\end{equation}
More precisely, they show that if $a(L)$ is a positive diverging function
of $L$, then
\begin{equation}
\P^{\frac{\ln(L)+a(L)}L}(\tilde\Theta\ge1)\to 0,
\qquad
\P^{\frac{\ln(L)-a(L)}L}(\tilde\Theta\ge1)\to 1.
\label{hegmar}
\end{equation}

In \cite{us}, we showed that there were of order of $L$ open minimal length
paths when the starting position $x$ is of order $1/L$ and we gave the
limiting law of
$\tilde\Theta/L$.

Informally, \eqref{hegmar} means 
that if the starting fitness~$x$ is larger than $(\ln L)/L$, then there
are no open minimal length path, and if $x$ is smaller than $(\ln L)/L$,
then there are some open minimal length paths. Even more informally,
the expectation \eqref{expminimal} tells the truth: when the expectation
goes to zero, there are no path (which is obvious); when the expectation
diverges, there are some paths (which is not automatic).

\subsection{Our results}

In this paper, we consider paths which are no longer of minimal length:
a mutation can change a 1 into a 0 as well as a 0 into a 1. We 
compute bounds for the expected number of open paths connecting $\sigO$ to
$\sigI$ given the starting fitness $x$ which lead to

\begin{theorem}
When $H=L$ (that is, when the fittest site is $\sigI$),
\begin{equation*}
\big[\E^x(\Theta)\big]^{1/L}\to \sinh(1-x) \qquad\text{as $L\to\infty$}.
\end{equation*}
In particular there is a critical value $x^*_1$ for the fitness of the
starting position,
\begin{equation*}
x^*_1=1-\arcsinh(1)=1-\ln(\sqrt2+1)=0.11863\ldots,
\end{equation*}
such that
\begin{itemize}
\item
For $x>x^*_1$, $\E^x(\Theta)$ goes to zero exponentially fast as $L\to\infty$
and, therefore $\P^x(\Theta\ge1)\to0$.
\item
For $x<x^*_1$, $\E^x(\Theta)$ diverges exponentially fast as $L\to\infty$.
\end{itemize}
As a consequence,
\begin{equation}
\limsup_{L\to\infty}\P(\Theta\ge1)\le x^*_1.
\end{equation}
\label{thm1}
\end{theorem}
We conjecture that the expectation ``tells the truth'' and that:
\begin{conjecture}
when $H=L$, for $x<x^*_1$,
\begin{equation*}
\P^x(\Theta\ge1)\to1 \qquad\text{as $L\to\infty$}.
\end{equation*}
and, as a consequence,
\begin{equation*}
\lim_{L\to\infty}\P(\Theta\ge1)=x^*_1.
\end{equation*}
\end{conjecture}

\begin{figure}[ht]
\centering\includegraphics[width=.7\textwidth]{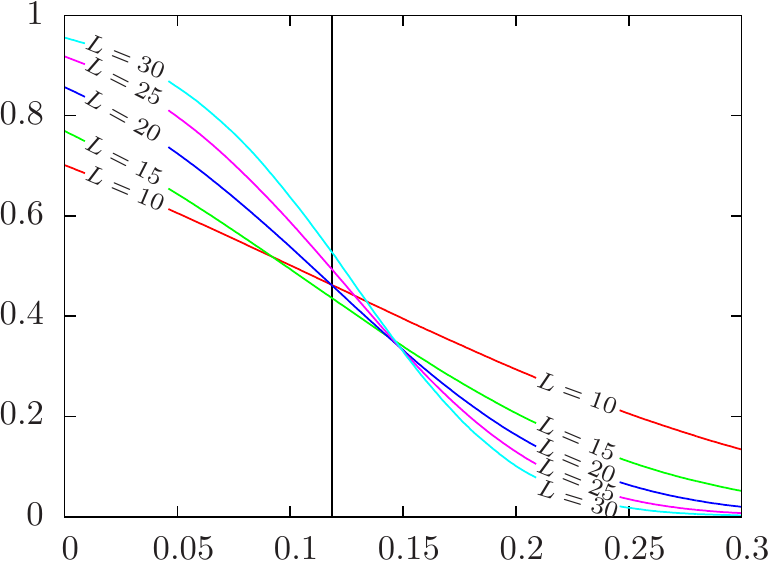}
\caption{Probability to have an open path from $\sigO$ to
$\sigI$ in the $L$-hypercube as a function of the starting
fitness~$x$. The critical $x^*_1$ is represented by the black
vertical line. These curves were obtained by Monte-Carlo
simulation on $10^6$ samples per size (only $10^5$ for the
largest size).}\label{figproba}
\end{figure}

As an illustration, Figure~\ref{figproba} shows the result of numerical
simulations measuring the probability  $\P^x(\Theta\ge1)$ that there are
some open paths on the $L$-hypercube as a function of $x$ for different
values of $L$. Our theorem is that, for large~$L$, the probability goes to
zero on the right of the black line and our conjecture is that it goes to
1 on the left. One might guess such a scenario from this picture alone,
with, however, a critical value around 0.15 rather than the actual
$x^*_1\approx0.12$.
Our work proves however that the critical value cannot be larger than
$x^*_1$.

When the fittest site is not $\sigI$, we have the following more general
result:
\begin{theorem}
Let $\alpha\in[0,1]$ and consider a function $L\mapsto H(L)$ such that
$H(L)/L\to\alpha$ as $L\to\infty$. Then, when $\sigend$ is chosen such
that its  Hamming distance is $H=H(L)$,
\begin{equation*}
\big[\E^x(\Theta)\big]^{1/L}\to \sinh(1-x)^\alpha\cosh(1-x)^{1-\alpha}
\qquad \text{as $L\to\infty$}.
\end{equation*}
In particular, for each $\alpha$, there is a critical value $x^*_\alpha$
for the fitness of the starting position, which is the unique solution of
\begin{equation*}
\sinh(1-x^*_\alpha)^\alpha \cosh(1-x^*_\alpha)^{1-\alpha}=1,
\end{equation*}
such that
\begin{itemize}
\item
For $x>x^*_\alpha$, $\E^x(\Theta)$ goes to zero exponentially fast as
$L\to\infty$ and, therefore $\P^x(\Theta\ge1)\to0$.
\item
For $x<x^*_\alpha$, $\E^x(\Theta)$ diverges exponentially fast as $L\to\infty$.
\end{itemize}
As a consequence,
\begin{equation}
\limsup_{L\to\infty}\P(\Theta\ge1)\le x^*_\alpha.
\end{equation}
\label{thm2}
\end{theorem}
We conjecture again that the expectation ``tells the truth'' and that:
\begin{conjecture}
when $H=H(L)$ with $H(L)/L\to\alpha$, for $x<x^*_\alpha$,
\begin{equation*}
\P^x(\Theta\ge1)\to1 \qquad\text{as $L\to\infty$}.
\end{equation*}
and, as a consequence,
\begin{equation*}
\lim_{L\to\infty}\P(\Theta\ge1)=x^*_\alpha.
\end{equation*}
\label{conj2}
\end{conjecture}

\begin{figure}[ht]
\centering\includegraphics[width=.6\textwidth]{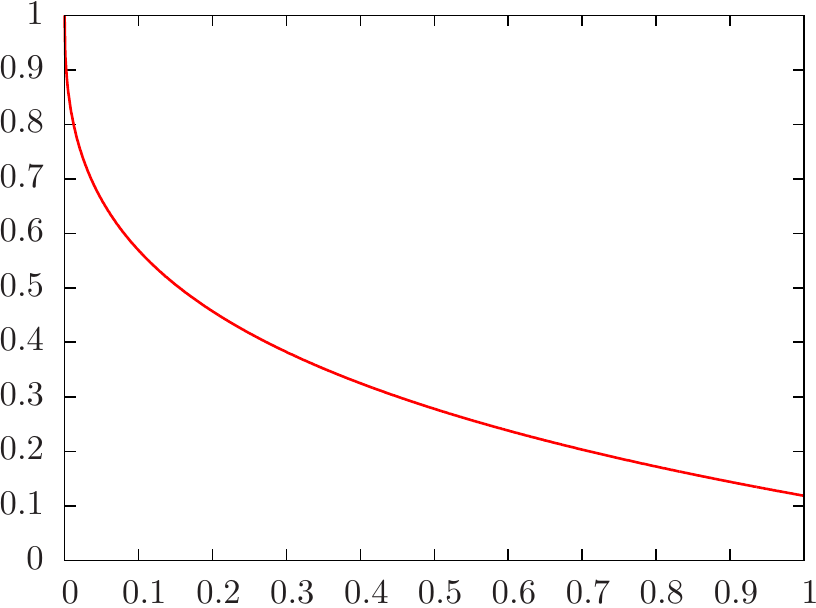}
\caption{The critical point $x_\alpha^*$ as a function of $\alpha$.}
\label{figcrit}
\end{figure}

Figure \ref{figcrit} gives the critical value $x_\alpha^*$ as a function of
$\alpha$. Noteworthy points are $x_1^*=1-\arcsinh(1)=0.118626\ldots$ as
already noted and
$x_{1/2}^*=1-\frac12\arcsinh(2)=1-\frac12\ln(2+\sqrt5)=0.278182\ldots$.

Finally, when the fittest site is chosen uniformly at random (this is the
model which is truly equivalent to the House-of-Cards model) it is clear
that for large $L$ the value of $H/L$ converges to $1/2$.
This leads to the following result:
\begin{theorem}
When the fittest site is chosen uniformly at random, one has
\begin{itemize}
\item For $x>x^*_{1/2}$, $\P^x(\Theta\ge1)\to0$ as $L\to\infty$.
\end{itemize}
As a consequence,
\begin{equation}
\lim_{L\to\infty}\P(\Theta\ge1)\le x^*_{1/2}.
\end{equation}
Furthermore, if one assumes Conjecture~\ref{conj2},
\begin{itemize}
\item For $x<x^*_{1/2}$, $\P^x(\Theta\ge1)\to1$ as $L\to\infty$; therefore
 $\P(\Theta\ge1)\to x^*_{1/2}$.
\end{itemize}
\label{thm3}
\end{theorem}

As a by-product of this work, we also found that the number of
self-avoiding paths (just plain paths, without any notion of fitness,
openness or accessibility) joining $\sigO$ to $\sigI$ in the $L$-hypercube
grows as a double exponential, see Theorem~\ref{thmaL} at the end of
Section~\ref{secaL}.

Theorem~\ref{thm1} if proved in Section~\ref{proof1}, Theorem~\ref{thmaL} in
Section~\ref{proofaL}, Theorem~\ref{thm2} in Section~\ref{later} and
Theorem~\ref{thm3} in Section~\ref{proof3}.

\section{Proof when the fittest site is $\sigI$}
\label{proof1}

We consider here the case $H=L$, \textit{i.e.}\@ when the fittest site, the one
with a fitness equal to 1, is $\sigI$. The generalization
to an arbitrary fittest site is described in Section~\ref{later}.

The minimum length of a path from $\sigO$ to $\sigI$ is $L$, as each one of
the $L$ bits has to be switched from 0 to 1. There exists however 
longer paths which have backsteps, \textit{i.e.}\@ steps where a bit is flipped
from 1 to 0. The length of a path with $p$ backsteps is clearly $L+2p$ as
each backstep must be compensated by an extra forward step.

We only need to consider paths that do not self-intersect, as it is obvious
that a path going twice to the same site cannot see its fitness increase
strictly.
We define
\begin{equation}\begin{aligned}
\all_{L}&=\text{the number of self-avoiding paths connecting
$\sigO$ to $\sigI$},\\
\all_{L,p}&=\text{the number of self-avoiding paths connecting
$\sigO$ to $\sigI$}\\&\quad\text{ with a length $L+2p$ (that is, with $p$ backsteps).}
\end{aligned}\end{equation}

\begin{figure}[ht]
\centering
\includegraphics[width=.7\textwidth]{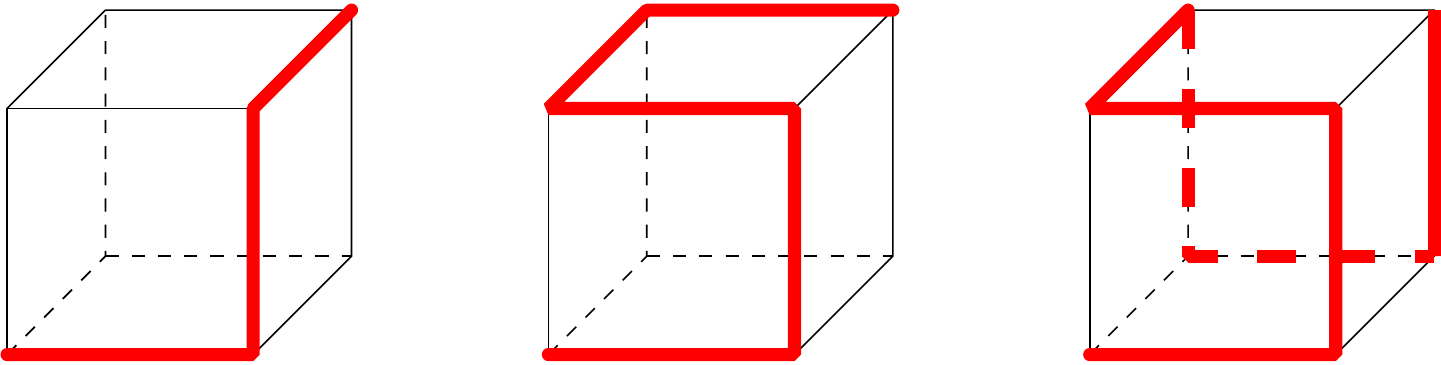}
\caption{The three self-avoiding paths on the cube connecting
$\sigO=$(bottom, left, front) to the opposite corner $\sigI=$(top, right,
back) which begin by
(right, up).}
\label{figcube}
\end{figure}

As an illustration, Figure~\ref{figcube} shows all self-avoiding paths
on the $3$-hypercube which begin by ``right, up'': there are three of them
with respective lengths 3 ($p=0$), 5 ($p=1$) and 7 ($p=2$). But there are
6 choices for the first two steps so, by symmetry, $\all_{3,0}=6$,
$\all_{3,1}=6$, $\all_{3,2}=6$ and, of course, $\all_3=18$.

When the starting site has a fixed fitness $x$, the probability that
a self-avoiding path of length $L+2p$ is open is simply (as in the
introduction)
$(1-x)^{L+2p-1}/(L+2p-1)!$: the $L+2p-1$ interior sites must be between $x$
and 1 and they must be in order.
Thus, the expected number of open paths when the starting value is $x$ is
\begin{equation}
\E^x(\Theta)=\sum_{p\ge0} \all_{L,p} \frac{(1-x)^{L+2p-1}}{(L+2p-1)!}.
\label{Ex}
\end{equation}

The problem therefore reduces to finding good estimates on $\all_{L,p}$.

\subsection{Some remarks about $\all_{L,p}$}
\label{secaL}

There seems to be little literature on the number $\all_{L}$ of
self-avoiding paths on a $L$-hypercube joining $\sigO$ to $\sigI$.
The sequence is referenced in the Online Encyclopedia
of Integer Sequences \cite{oeis} and values are given up to $L=5$:
\begin{equation}\begin{gathered}
\all_1=1,\quad
\all_2=2,\quad
\all_3=3!\times3=18,\quad\\
\all_4=4!\times268=6\,432,\quad
\all_5=5!\times155\,429\,607=18\,651\,552\,840.
\end{gathered}
\end{equation}
Furthermore, the value for $\all_5$ is mentioned in an evolutionary
biology paper \cite{DePristo.2007}. These numbers were obtained by brute
force computer enumeration; due to the combinatorial explosion, $\all_6$
is out of reach by this method.

For a given $L$, the length of a self-avoiding path cannot exceed the
number $2^L-1$ of sites to explore; hence one must have $L+2p\le 2^L-1$.
Just to give a flavour of the structure of paths, here are all the numbers
$\all_{5,p}$:
\begin{equation}\begin{gathered}
\all_{5,0}=5!,\quad
\all_{5,1}=5!\times10,\quad
\all_{5,2}=5!\times107,\quad
\all_{5,3}=5!\times1\,097,\quad
\all_{5,4}=5!\times9\,754,\\
\all_{5,5}=5!\times72\,305,\quad
\all_{5,6}=5!\times448\,536,\quad
\all_{5,7}=5!\times2\,243\,671,\\
\all_{5,8}=5!\times8\,631\,118,\quad
\all_{5,9}=5!\times24\,044\,702,\quad
\all_{5,10}=5!\times44\,617\,008,\\
\all_{5,11}=5!\times48\,280\,086,\quad
\all_{5,12}=5!\times24\,000\,420,\quad
\all_{5,13}=5!\times3\,080\,792.
\end{gathered}
\end{equation}

It is clear that $\all_{L,p}$ must be a multiple of $L!$ as from
a given path more paths can be built by simply applying a permutation of
the $L$ directions of the edges. For $p=0$, one has obviously
$\all_{L,0}=L!$. 
For $p=1$, the paths have length $L+2$ with one
backstep at some position $k$ with $3\le k\le L$. (The first step cannot be
a backstep as all the bits are still 0. The second step cannot be
a backstep as it would bring the path back to the starting position. The
steps
$L+1$ and $L+2$ cannot be the backstep as the system would have already
reached the end point at step $L$.) For a given backstep position $k$,
there are $L!/(L-k+1)!$ possible paths up to step $k-1$, and $k-2$
possible choice for the backstep (because there are $k-1$ bits set to 1 but
one cannot choose the bit that was just set), and $L-k+1$ choices for
the step following the backstep (because there are $L-k+2$ bits set to
0 but one cannot choose the bit that was just unset) and $(L-k+1)!$ choices
for all the subsequent steps. The number of paths with a given backstep
position $k$ is then $L!\times(k-2)(L-k+1)$; summing over all $k$ gives
\begin{equation}
\all_{L,1}=L!\times\frac{L(L-1)(L-2)}6.
\end{equation}

A similar (but much more strenuous) derivation leads for $p=2$ to
\begin{equation}
\all_{L,2}=
L! \times\frac{ (L - 1)(L - 2)(5 L^4 + 3 L^3 + 34 L^2 - 264 L + 180)}{360}.
\label{allL2}
\end{equation}

It is easy to convince oneself that
for fixed $p$, as $L\to\infty$,
\begin{equation}
\all_{L,p}\sim L! \frac{L^{3p}}{6^pp!}.
\label{asymp1}
\end{equation}
Indeed, one needs to choose $p$ backsteps at positions
$3\le k_1<k_2<\cdots<k_p\le L+2p-2$
in a sequence of $L+2p\approx L$ steps (we are dropping all the
non-dominant terms). At the $j$-th backstep there are of order $k_j$
choices to choose the bit we set to 0 (actually: $k_j-2j$ choices if the
previous step was not a backstep, but $k_j\propto L$ and $j\le p$ and we
are dropping all the
non-dominant terms). The step after backstep~$j$ has to leading order
$L-k_j$ bits~0 which can be switched to~1, and all the other steps combine to
build $L!$. With the $k_j$ given, one therefore gets a number of
paths of order $L!\times k_1(L-k_1)\times k_2(L-k_2)\times\cdots\times
k_p(L-k_p)$. Summing over the ordered $k_j$'s leads to \eqref{asymp1}. The
expression for fixed $k_j$'s is far from being correct if there are
backsteps too close to each other but all of this only contribute to
the next order term in the expression.

It would be extremely interesting to understand better how
$\all_L$ grows with $L$.
In this work, we give upper and lower
bounds for $\all_{L,p}$ to study the problem defined
in the introduction which, as a by-product, also lead to the following
theorem proved in Section~\ref{proofaL}.
\begin{theorem}\label{thmaL}
Recall that $\all_L$ is the number of self-avoiding paths on the
$L$-hypercube from $\sigO$ to $\sigI$. Then
\begin{equation*}
\lim_{L\to\infty}\frac{\ln\ln \all_L}{L}=\ln 2.
\end{equation*}
More precisely, there exists two positive constants $c$ and $c'$ such
that, for $L$ large enough,
\begin{equation*}
c \le \frac{\ln \all_L }{2^L} \le c'\ln L.
\end{equation*}
\end{theorem}
 
\subsection{Coding of a path}
\label{coding}

We use the following representation for a path on the hypercube:
it is a string of numbers between $1$ and $L$ where each
number indicates the position of the bit being flipped by the corresponding
step. The first time a particular number is met, the bit is flipped from
0 to 1; the next time from 1 to 0, etc. To take an example, the paths of
Figure~\ref{figcube} would be respectively coded ``123'', ``12131'' and
``1213212'' assuming that 1 is the left/right direction, 2 is the up/down
direction and 3 is the front/back direction.

Clearly, for paths going from $\sigO$ to $\sigI$, each number must appear
an odd number of times in the string. A path visits twice the same site if
there exists a non-empty substring\footnote{We recall that a substring is
a subsequence of consecutive terms.} of the path where each number appears an
even number of times (including zero, of course), as all the bits encoding
the position are clearly the same before and after the substring.
Examples of minimal forbidden substrings include ``11'', ``1212'',
``12313424'', etc. A self-avoiding path is then, of course, such that
there is no such substring. 

\subsection{Upper bound}
\label{upper}
In this section we prove the following:
\begin{lemma}
When $H=L$,
\begin{equation*}
\E^x(\Theta)\le L \sinh(1-x)^L\cotanh(1-x).
\end{equation*}
\label{lemupper}
\end{lemma}
This implies that $\limsup_{L\to\infty}
\big[\E^x(\Theta)\big]^{1/L}\le\sinh(1-x)$, which is the first half of the
proof of Theorem~\ref{thm1}.

Let
\begin{equation}\begin{aligned}
\majo_{L,p}&=\text{the set of paths connecting
$\sigO$ to $\sigI$ with a}
\\&\quad\text{ length $L+2p$ (that is, with $p$ backsteps)
 where intersections are allowed.}
\end{aligned}
\end{equation}

Clearly, 
\begin{equation}
\all_{L,p}\le \majo_{L,p},
\label{allmajo}
\end{equation}
(by an abuse of notation, the cardinal of $\majo_{L,p}$ is also noted
$\majo_{L,p}$) and it turns out that this very simple
upper bound is sufficient for our purpose.

With our representation, $\majo_{L,p}$ is the set of strings with length
$L+2p$ made of numbers between $1$ and $L$ where each number appears an odd
number of times. 
We build $\majo_{L,p}$ by recurrence: for any $p$,
there is one path in $\majo_{1,p}$: it is the path ``111\ldots'' where one
walks back and forth between the two sites of the $1$-hypercube.

To construct a path in $\majo_{L+1,p}$ (with length $L+1+2p$), we
\begin{itemize}
\item choose how many times the number $L+1$ appears. This number is odd,
let it be $2q+1$ with $0\le q\le p$,
\item choose the positions of the $2q+1$ numbers $L+1$ amongst the $L+1+2p$
possible positions in the string,
\item fill in the remaining $L+2p-2q$ positions with the string coding an
arbitrary path chosen in $\majo_{L,p-q}$.
\end{itemize}

In equations, this construction gives
\begin{gather}
\majo_{1,p}=1,\qquad
\majo_{L+1,p}=\sum_{q=0}^p\binom{L+1+2p}{2q+1} \majo_{L,p-q}.
\label{recurmajo}
\end{gather}
Writing the binomial with factorials one gets
\begin{equation}
\frac{\majo_{L+1,p}}{(L+1+2p)!}=\sum_{q=0}^p
\frac{\majo_{L,p-q}}{(L+2(p-q))!}\times\frac1{(2q+1)!}.
\end{equation}

Let $G_L(X)$ be the generating function defined by
\begin{equation}
G_L(X):=\sum_{p\ge0}\frac{\majo_{L,p}}{(L+2p)!} X^{L+2p}.
\label{genmajo}
\end{equation}
Notice that from \eqref{Ex} and \eqref{allmajo}, one has
\begin{equation}
\E^x(\Theta)\le G_L'(1-x).
\label{upb}
\end{equation}

The recurrence on $M_{L,p}$ translates into a recurrence on $G_L(X)$:
\begin{equation}
\begin{aligned}
G_{L+1}(X)
&=
\sum_{p\ge0}\sum_{q=0}^p\frac{\majo_{L,p-q}}{(L+2(p-q))!}X^{L+2(p-q)}\times
\frac{X^{2q+1}}{(2q+1)!},
\\&=
\sum_{q\ge0}\sum_{p\ge q}
\frac{\majo_{L,p-q}}{(L+2(p-q))!}X^{L+2(p-q)}\times
\frac{X^{2q+1}}{(2q+1)!},
\\&=
\sum_{q\ge0}\sum_{p\ge 0}
\frac{\majo_{L,p}}{(L+2p)!}X^{L+2p}\times
\frac{X^{2q+1}}{(2q+1)!},
\\&=
G_L(X)\sinh(X).
\end{aligned}
\label{computeGL}
\end{equation}
But $G_1(X)=\sinh(X) $, hence 
\begin{equation}
G_L(X)=\sinh(X)^L
\label{G_L}
\end{equation}
which, with \eqref{upb}, concludes the proof of Lemma~\ref{lemupper}.

\subsection{Lower bound}
\label{lower}
In this section we proove the second half of Theorem~\ref{thm1}:
\begin{lemma}
When $H=L$,
\begin{equation*}
\limsup_{L\to\infty} \big[\E^x(\Theta)\big]^{1/L}\ge\sinh(1-x).
\end{equation*}
\label{lemlower}
\end{lemma}

For this lower bound, we construct a subset  $\mino_{L}$ of all the
self-avoiding paths on the $L$-hypercube joining $\sigO$ to $\sigI$. The
construction is recursive:
\begin{itemize}
\item For $L=1$, there is only one self-avoiding path joining the two
corners of the $1$-hypercube. In our coding, this path is represented by
the string ``1''.
\item A path (coded by a string) is in $\mino_{L+1}$ if\quad (a)~the number
$L+1$ appears an odd number of times in the string,\quad(b)~the number
$L+1$ never appears at two consecutive positions, and\quad(c)~the string
where all the number $L+1$ are removed codes a path in $\mino_L$.
\end{itemize}
It is clear by recurrence that paths in $\mino_L$ are valid self-avoiding
paths on the $L$-hypercube.
%
To illustrate, $\mino_1=\{``1"\}$, $\mino_2=\{``12"$, $``21"\}$, $\mino_3=\{``123"$, 
$``132"$, $``312"$, $``213"$, $``231"$, $``321"$, $``31323"$, $``32132"\}$.
Notice how we lost symmetry in the paths: $``31323"$ is in $\mino_3$ but not
$``13121"$.

We now define
\begin{equation}
\mino_{L,p}=\text{the set of paths of length $L+2p$ in $\mino_L$}.
\end{equation}
Clearly
\begin{equation}
\mino_{L,p}\le \all_{L,p},
\end{equation}
(once again, by an abuse of notation we write the cardinal of $\mino_{L,p}$
also as $\mino_{L,p}$). 

One has $\mino_{1,0}=\{``1"\}$ and $\mino_{1,p}=\emptyset$ for $p>0$. 
The sets $\mino_{L,p}$ are then built recursively:  to construct a path in
$\mino_{L+1,p}$, we
\begin{itemize}
\item choose how many times the number $L+1$ appears. This number is odd,
let it be $2q+1$ with $0\le q\le p$,
\item choose the positions of the $2q+1$ numbers $L+1$ amongst the $L+1+2p$
possible positions in the string in such a way that there are no two
consecutive numbers $L+1$,
\item fill in the remaining $L+2(p-q)$ positions with the string coding
an arbitrary path chosen in $\mino_{L,p-q}$.
\end{itemize}

Let us recall that the number of ways of choosing $P$ items out of
a sequence of $N$ such that two consecutive items in the sequence
cannot be both chosen is $\binom{N-P+1}P$. Indeed, each configuration
can be bijectively obtained by first choosing $P$ items out of a sequence
of $N-P+1$ and then expanding the sequence to size $N$ by inserting
one unchosen item ``$\cdot$'' before each chosen item~``$\bullet$'' except
the first one; for instance, with $P=4$ and $N=11$: ($\bullet\cdot
\bullet \bullet\cdot\cdot \bullet \cdot\to \bullet\cdot\cdot \bullet\cdot
\bullet\cdot\cdot\cdot \bullet\cdot$).

With this result, the construction of $\mino_{L,p}$ leads to
\begin{equation}
\mino_{1,p}=\one_{p=0},\qquad \mino_{L+1,p}=\sum_{q=0}^p
\binom{L+1+2p-2q}{2q+1}\mino_{L,p-q}.
\label{recurmino}
\end{equation}

In order to write a generating function similar to the $G_L(X)$ defined in
the previous section, we need to replace the binomial in \eqref{recurmino} by
the same binomial as in \eqref{recurmajo}. Let us write, for $0\le q\le p$,
\begin{equation}
\begin{aligned}
\binom{L+1+2p-2q}{2q+1}
&=\frac1{(2q+1)!}\times\quad(L+1+2p-2q)!\quad\times\frac1{(L+2p-4q)!}
\\
&=\frac1{(2q+1)!}\times\frac{(L+1+2p)!}{\prod_{k=0}^{2q-1} (L+1+2p-k)}
\times\frac{\prod_{k=0}^{2q-1} (L+2p-2q-k)}{(L+2p-2q)!}
\\
&=\binom{L+1+2p}{2q+1}\times
\prod_{k=0}^{2q-1}\frac{L+2p-2q-k}{L+1+2p-k}
\\&=\binom{L+1+2p}{2q+1}\times\prod_{k=0}^{2q-1}
\Big[1-\frac{2q+1}{L+1+2p-k}\Big]
\\&\ge\binom{L+1+2p}{2q+1}\times
\Big[1-\frac{2q+1}{L+2}\Big]^{2q}\one_{2q<L+1}\qquad\text{(using $p\ge q$).}
\end{aligned}
\end{equation}

Then, defining $\tilde \mino_{L,p}$ by the recurrence
\begin{equation}
\tilde \mino_{1,p}=\mino_{1,p}=\one_{p=0},\qquad
\tilde \mino_{L+1,p}=\sum_{q=0}^p
\tilde \mino_{L,p-q}\binom{L+1+2p}{2q+1}\times \Big[1-\frac{2q+1}{L+2}
\Big]^{2q}
\one_{2q< L+1}
\end{equation}
it is clear that 
\begin{equation}
\tilde \mino_{L,p}\le \mino_{L,p}\le \all_{L,p}.
\label{allmino}
\end{equation}

As for the upper bound, let $g_L(X)$ be the generating funcion of the
$\tilde\mino_{L,p}$ defined by the finite sum
\begin{equation}
g_L(X):=\sum_{p\ge0} \frac{\tilde \mino_{L,p}}{(L+2p)!}X^{L+2p}.
\label{genmino}
\end{equation}
Notice that from \eqref{Ex} and \eqref{allmino}, one has
\begin{equation}
g_L'(1-x) \le\E^x(\Theta).
\label{lob}
\end{equation}
The recurrence on $\tilde\mino_{L,p}$ translates into a recurrence on
$g_L(X)$. One gets easily, by the same argument  as in \eqref{computeGL}
\begin{equation}
g_{L+1}(X)=g_{L}(X)\times\sum_{q\ge0} \frac{X^{2q+1}}{(2q+1)!}
\Big[1-\frac{2q+1}{L+2}\Big]^{2q}\one_{2q<L+1}.
\label{reca}
\end{equation}
Defining
\begin{equation}
\sinh_l(X):=\sum_{q\ge0}
\frac{X^{2q+1}}{(2q+1)!} \Big[1-\frac{2q+1}{l+1}\Big]^{2q}
\one_{2q<l},
\label{33}
\end{equation}
then
\eqref{reca} reads $g_{L+1}(X)=g_L(X)\sinh_{L+1}(X)$. Furthermore,
$g_1(X)=X=\sinh_1(X)$ so that
\begin{equation}
g_L(X)=\prod_{l=1}^L \sinh_l(X).
\end{equation}
The derivative of $g_L$ can then be written
\begin{equation}
g_L'(X)=
\Big(\sum_{l=1}^L\frac{\sinh_l'(X)}{\sinh_l(X)}\Big)\times g_L(X)
\label{deriv}
\end{equation}
It is clear that $X\le\sinh_l(X)\le\sinh(X)$ and that
$1\le\sinh_l'(X)\le\cosh(X)$. The sum in \eqref{deriv} is therefore bounded
between $L/\sinh(X)$ and $L\cosh(X)/X$, and the sum to the power $1/L$
converges to 1 as $L\to\infty$.

Furthermore, by dominated convergence,
$\sinh_l(X)\to\sinh(X)$ (and $\sinh_l'(X)\to\cosh(X)$) as $l\to\infty$.
This is sufficient to imply that $g_L(X)^{1/L}$ converges by Cesaro to
$\sinh(X)$. Finally, $g_L'(X)^{1/L}\to\sinh(X)$ and, with~\eqref{lob},
\begin{equation}
\sinh(1-x)\le\liminf_{L\to\infty}\big[\E^x(\Theta)\big]^{1/L},
\end{equation}
which is the second half of Theorem~\ref{thm1}.

\section{Proof of Theorem~\ref{thmaL}}
\label{proofaL}

In the previous section, we used the bounds $\mino_{L,p}\le \all_{L,p}\le
\majo_{L,p}$
to obtain an estimate on $\E^x(\Theta)$ with \eqref{Ex}. We now use the same
bounds to obtain an estimate on $\all_L=\sum_p \all_{L,p}$ and prove
Theorem~\ref{thmaL}.

\paragraph{lower bound}

We define another generating function of the $\mino_{L,p}$. Let
\begin{equation}
\phi_L(X)=\sum_{p\ge0}\mino_{L,p}X^{L+2p}.
\end{equation} 
Then one finds easily from \eqref{recurmino} that
\begin{equation}
\phi_1(X)=X,\qquad
\phi_{L+1}(X)=\frac{1+X}2\phi_L\big(X+X^2\big)
             -\frac{1-X}2\phi_L\big(X-X^2\big).
\label{recphi}
\end{equation}
From its definition $\phi_L(X)$ is a polynomial in $X$ (recall that
$\all_{L,p}$ and, therefore, $\mino_{L,p}$ is zero
if $p$ is too large). This polynomial is an odd function of~$X$
if $L$ is odd and an even function if $L$
is even. Let $d_L$ be the degree of this polynomial. By considering the
highest degree term in \eqref{recphi} one gets easily
\begin{equation}
d_1=1,\qquad
d_{L+1}=\begin{cases}
2d_L&\text{if $L$ is odd,}\\
2d_L+1&\text{if $L$ is even.}
\end{cases}
\end{equation}
This can be solved into
\begin{equation}
d_L=
\begin{cases}
\frac{2^{L+1}-1}{3}&\text{if $L$ is odd,}\\[1ex]
\frac{2^{L+1}-2}{3}&\text{if $L$ is even.}
\end{cases}
\end{equation}

By definition, $\mino_L=\sum_p \mino_{L,p}=\phi_L(1)$. From \eqref{recphi}
one has, furthermore,
\begin{equation}
\mino_L=\phi_L(1)=\phi_{L-1}(2).
\end{equation}
As the polynomials $\phi_L$ have non-negative integer coefficients, one
clearly have
\begin{equation}
\all_L\ge \mino_L\ge 2^{d_{L-1}}
\end{equation}
which is enough for the lower bound of Theorem~\ref{thmaL}.

\paragraph{Upper bound}

There are infinitely many paths in $\majo_L$, but one knows that paths
in $\all_L$ have a maximum length of $2^L$ (one could be more precise: at most
 $2^L-1$ if $L$ is odd and at most $2^L-2$ if $L$ is even) so that
\begin{equation}
\all_L\le\sum_p\majo_{L,p}\one_{L+2p\le 2^L}.
\end{equation}

Let us write for an analytical function $f$
\begin{equation}
T_nf=\text{The Taylor polynomial of $f$ of degree $n$.}
\end{equation}
Our upper bound is then
\begin{equation}
\all_L\le\int_0^\infty e^{-X} T_{2^L}G_L(X)\,\diffd X,
\label{44}
\end{equation}
as seen from the definition \eqref{genmajo} of $G_L$.

For any absolutely increasing function~$f$ (all derivatives are
non-negative), any order $n$ and any cutoff point $C$, one has
\begin{equation}
T_nf(X)\le\begin{cases}
f(X)&\text{if $X\le C$},\\
f(C)\left(\frac{X}C\right)^n&\text{if $X\ge C$}.
\end{cases}\end{equation}
The first line is trivial as we removed some non-negative terms. The second
line is also trivial because $X^k\le C^k (X/C)^n$ for $X\ge C$ and $k\le
n$, so the inequality holds for all the terms in the polynomial $T_nf$.
Then
\begin{equation}
\int_0^\infty e^{-X} T_nf(X)\,\diffd X\le f(C) C+ f(C)\frac{n!}{C^n}.
\end{equation} 

We apply this to $f(X)=G_L(X)=\sinh(X)^L\le e^{LX}/2^L$ and to $C=n/L$:
\begin{equation}
\int_0^\infty e^{-X} T_nG_L(X)\,\diffd X
\le \frac{n e^n}{L 2^L} + \frac{n! e^n}{n^n}\frac{L^n}{2^L}.
\end{equation}
Remember that by Stirling $e^n n!/n^n\sim\sqrt{2\pi n}$. The second term
on the right hand-side
is much larger than the first ($L^n$ vs $e^n$). Replace $n$ by
$2^L$, and it is easy to check that the bound~\eqref{44} gives the second
half of Theorem~\ref{thmaL}.

\section{Proof for an arbitrary $H$}
\label{later}

We assume now that the fittest site, the one with a fitness equal to 1, is
no longer $\sigI$ but rather an arbitrary given site~$\sigend$. By
symmetry, the accessibility of $\sigend$ depends only on the number of bits
set to~1 in $\sigend$; let $H$ (as in ``Hamming distance'') be this number
of bits. To simplify the discussion, we assume that the bits $1$ to $H$ in
$\sigend$ are set to~1 and that the bits $H+1$ to $L$ are set to~0.

We emphasize that we consider paths on the $L$-hypercube and not on the
$H$-hypercube: valid paths may leave the $H$-hypercube and then do
backsteps to go back to $\sigend$. In previous studies where only shortest
paths where considered (no backstep), considering the case
$\sigend\ne\sigI$ was meaningless as it was 
simply equivalent to changing the dimension of the hypercube.

The minimum length of a path from $\sigO$ to $\sigend$ is $H$. A path with
$p$ backsteps has a length $H+2p$. We define
\begin{equation}\begin{aligned}
\all_{L,H,p}&=\text{the number of self-avoiding paths connecting
$\sigO$ to $\sigend$}\\&\quad\text{ with a length $H+2p$ (that is, with $p$
backsteps).}
\end{aligned}\end{equation}
Then, with the same argument as before,
\begin{equation}
\E^x(\Theta)=\sum_{p\ge0}\all_{L,H,p}\frac{(1-x)^{H+2p-1}}{(H+2p-1)!}.
\end{equation}

With the coding introduced in Section~\ref{coding}, a self-avoiding
path from $\sigO$ to $\sigend$ is a string of numbers between $1$ and $L$
such that\quad(a)~the numbers between 1 and~$H$ appear an odd number of
times,\quad(b)~the numbers between~$H+1$ and $L$ appear an even number of
times (including zero),\quad(c)~in any non-empty substring, there must be
at least one number which appears an odd number of times.

Using the same strategy as in the previous section, we bound
$\all_{L,H,p}$:
\begin{equation}
\mino_{L,H,p}\le\all_{L,H,p}\le\majo_{L,H,p},
\end{equation}
where
\begin{itemize}
\item
$\majo_{L,H,p}$ is the number (or the set) of paths on the $L$-hypercube
of length $H+2p$ from
$\sigO$ to $\sigend$  where intersections are authorized. Naturally,
$\majo_{H,H,p}=\majo_{H,p}$ (as defined in Section~\ref{upper}).
For $L\ge H$, we obtain recursively $\majo_{L+1,H,p}$ as the number of
strings of length $H+2p$ such that $L+1$ appears an even number of times
$2q$ 
\emph{at arbitrary positions} and such that if one removes all the
occurrences of $L+1$, the resulting string is in $\majo_{L,H,p-q}$.
\item $\mino_{L,H,p}$ is defined recursively: for $H=L$ one has
$\mino_{H,H,p}=\mino_{H,p}$
(as defined in Section~\ref{lower}). For $L\ge H$, we obtain
recursively $\mino_{L+1,H,p}$ as
the number of strings of length $H+2p$ such that $L+1$ appears an even
number of times $2q$ \emph{but never at two consecutive positions} and such
that if one removes all the occurrences of $L+1$, the resulting string is in
$\mino_{L,H,p-q}$.
\end{itemize}

These definitions translate directly into the following equations for
$L\ge H$:
\begin{equation}
\majo_{L+1,H,p}=\sum_{q=0}^p \binom{H+2p}{2q} \majo_{L,H,p-q},
\qquad
\mino_{L+1,H,p}=\sum_{q=0}^p \binom{H+2p-2q+1}{2q} \mino_{L,H,p-q},
\end{equation}
to be compared with \eqref{recurmajo} and \eqref{recurmino}.
In each case, $2q$ is
the number of times $L+1$ appears in the string of length $H+2p$. The two
binomials correspond respectively to the number of ways of choosing $2q$
elements in $H+2p$, and the number of ways of choosing $2q$ in $H+2p$
such that two consecutive elements may not be chosen.

For the upper bound we define as before the generating function
$G_{L,H}(X)$:
\begin{equation}
G_{L,H}(X):=\sum_{p\ge0} \frac{\majo_{L,H,p}}{(H+2p)!} X^{H+2p}.
\end{equation}
Using the same technique as in \eqref{computeGL}, one gets, for $L\ge H$,
\begin{equation}
G_{L+1,H}(X)=G_{L,H}(X)\cosh(X).
\end{equation}
Furthermore, since  $G_{H,H}(X)=G_H(X)=\sinh(X)^H$, one has
\begin{equation}
G_{L,H}(X)=\sinh(X)^H \cosh(X)^{L-H}.
\label{55}
\end{equation}

Now, for the lower bound, we make the same transformation as before and
we obtain
\begin{equation}
\binom{H+2p-2q+1}{2q}
\ge
\binom{H+2p}{2q} \Big[1-\frac{2q}{H+1}\Big]^{2q-1} \one_{2q< H+1}.
\end{equation}
Obviously by recurrence, $\mino_{L,H,p}\ge\tilde\mino_{L,H,p}$ where we
define
\begin{equation}
\tilde\mino_{H,H,p}=\tilde\mino_{H,p},\qquad
\tilde\mino_{L+1,H,p}=\sum_{q=0}^p\binom{H+2p}{2q}\tilde\mino_{L,H,p-q} 
\Big[1-\frac{2q}{H+1}\Big]^{2q-1}
\one_{2q< H+1}.
\label{recurminoH}
\end{equation}
Introducing the generating function $g_{L,H}(X)$
\begin{equation}
g_{L,H}(X):=\sum_{p\ge0}\frac{\tilde\mino_{L,H,p}}{(H+2p)!}X^{H+2p},
\end{equation}
one gets from \eqref{recurminoH}
\begin{equation}
g_{H,H}(X)=g_H(X)=\prod_{l=1}^H\sinh_l(X),\qquad
g_{L+1,H}(X)=g_{L,H}(X)\cosh_H(X),
\end{equation}
where $\sinh_l(X)$ was defined in \eqref{33} and where
\begin{equation}
\cosh_H(X)=\sum_{q\ge0}\frac{X^{2q}}{(2q)!}
\Big[1-\frac{2q}{H+1}\Big]^{2q-1}\one_{2q< H+1}.
\end{equation}
This gives
\begin{equation}
g_{L,H}(X)=\Big(\prod_{l=1}^H\sinh_l(X)\Big)\cosh_H(X)^{L-H}.
\end{equation}
Finally, collecting the bits, one has the bounds
\begin{equation}
g_{L,H}'(1-x)\le\E^x(\Theta)\le G_{L,H}'(1-x).
\label{62}
\end{equation}

If one chooses a function $L\mapsto H(L)$ such that  $H(L)/L\to\alpha$ as
$L$ goes to infinity, then it is very easy to check that
$\big[g'_{L,H(L)}(1-x)\big]^{1/L}$ and
$\big[G'_{L,H(L)}(1-x)\big]^{1/L}$ both converges as $L\to\infty$ to the
same quantity $\sinh(1-x)^\alpha\cosh(1-x)^{1-\alpha}$, which proves
Theorem~\ref{thm2}.
\section{Proof when the fittest point is random}
\label{proof3}

In this short section we prove Theorem \ref{thm3}.
We consider the situation in which  the fittest site $\sigend$  is
chosen uniformly at random in the hypercube (this is the case which  is
equivalent to the House-of-Cards model). In this case,
$H$ is obviously a binomial with parameters $L$ and $1/2$, and therefore,
$\alpha := H/L \to \frac12$ in probability.
Thus,  when $x >x^*_{1/2},$ for $\epsilon >0$ small enough so that
$x>x^*_{1/2-\epsilon}$ we have that
\begin{equation}
\begin{aligned}
\P^x(\Theta \ge 1)  &\le \P^x ( \Theta\ge 1 , \alpha \in \frac12 \pm
\epsilon) +\P^x(\alpha \not \in \frac12 \pm \epsilon),  \\
& \to 0\qquad\text{as $L\to\infty$}.
\end{aligned}\label{75}
\end{equation}
To see this, just observe that the second term on the right-hand side tends
to 0 independently of the value of $x$ while the first term also goes to
0 as for any $\alpha >1/2 -\epsilon$ we always have $x > x^*_{\alpha}$.

Moreover, if we assume that Conjecture \ref{conj2} holds, \textit{i.e.}\@ that for all
$\alpha$ fixed, $\P^x(\Theta \ge 1) \to 1$ when $x <x^*_{\alpha}$, the
second part of Theorem \ref{thm3} follows by the same argument.

However, be wary that the expected number of paths is lying. Just
looking at the upper bound (but the lower bound should be the same) one
can write with \eqref{55} and~\eqref{62}:
\begin{equation}
\E^x(\Theta)\le\sum_{H=0}^L\frac{1}{2^L}\binom{L}{H}G_{L,H}'(1-x)=L
\Big[\frac{e^{1-x}}2\Big]^L,
\end{equation}
which  diverges exponentially if and only if
$x\le1-\ln2=0.30685\ldots$. This seems to
give a critical point which is larger than $x^*_{1/2}$, but 
what happens is that with an
exponentially small probability, $\alpha$ is much smaller than $1/2$ which
generates exponentially many paths thereby contributing to the expectation.

\end{document}